\documentclass[11pt]{amsart}

\usepackage[T1]{fontenc}
\usepackage{times}
\usepackage{amssymb,epsfig,verbatim,xypic}
\usepackage{calligra,mathrsfs}
\theoremstyle{plain}
\usepackage[T1]{fontenc}
\usepackage[utf8]{inputenc}
\usepackage[english]{babel}
\usepackage{amssymb}
\usepackage{amsthm}
\usepackage{graphics}
\usepackage{amsmath}
\usepackage{amstext}
\usepackage{multirow}
\usepackage{enumerate}

\usepackage{hyperref}
\usepackage{graphicx}
\usepackage{amsmath}
\usepackage{amssymb}
\usepackage{amsthm}
\usepackage{amstext}
\usepackage{epstopdf}
\usepackage{mathrsfs}
\usepackage{amscd}
\usepackage{tikz}
\usetikzlibrary{cd}
\AtBeginDocument{%
   \def\MR#1{}
}

\newtheorem{thm}{Theorem}[section]
\newtheorem{lemma}[thm]{Lemma}
\newtheorem{prop}[thm]{Proposition}

\newtheorem{cor}[thm]{Corollary}
\newtheorem{question}[thm]{Question}

\newtheorem{THM}{Theorem}

\newtheorem{ex}[thm]{Example}

\theoremstyle{remark}
\newtheorem{remark}[thm]{Remark}

\newcommand{\mb}{\mathbb}

\newcommand{\mc}{\mathcal}

\newcommand{\C}{\mb C}

\newcommand{\F}{\mc F}
\newcommand{\G}{\mc G}

\newcommand\restr[2]{{
  \left.\kern-\nulldelimiterspace 
  #1 
  \vphantom{\big|} 
  \right|_{#2} 
  }}

\DeclareMathOperator{\codim}{codim}
\DeclareMathOperator{\sing}{sing}

\DeclareMathOperator{\Aut}{Aut}

\DeclareMathOperator{\trdeg}{tr\, deg_{\mathbb C}}

\DeclareMathOperator{\HomSheaf}{\mathscr{H}\text{\kern -3pt {\calligra\large om}}\,}

\newcommand{\ie}{{\it{i.e. }}}

\numberwithin{equation}{section}

\addtocounter{section}{0}             
\numberwithin{equation}{section}       

\title[Submanifolds with ample normal bundle]{Submanifolds with ample normal bundle}

\author{Maycol Falla Luza}
\address{Universidade Federal Fluminense, Rua Prof. Marcos Waldemar de Freitas Reis, S/N,  Niterói,  24210-201, Rio de Janeiro, Brazil}

\author{Frank Loray}
\address{Univ Rennes, CNRS, IRMAR - UMR 6625, F-35000 Rennes, France.}

\author{Jorge Vitório Pereira}
\address{IMPA, Estrada Dona Castorina, 110, Horto, 22460-320,  Rio de
Janeiro, Brasil}

\date{\today}

\begin{document}

\begin{abstract}
    We construct germs of complex manifolds of dimension $m$ along projective submanifolds of dimension $n$ with ample normal bundle and without non-constant meromorphic functions whenever $m \ge 2n$. We also show that our methods do not allow the construction of similar examples when $m < 2n$  by establishing an algebraicity criterion for foliations on projective spaces which generalizes a classical result by Van den Ven characterizing linear subspaces of projective spaces as the only submanifolds with split tangent sequence.
\end{abstract}

\maketitle
\setcounter{tocdepth}{1}
\tableofcontents

\section{Introduction}

Let $X$ be a (not necessarily compact) complex manifold and let $Y\subset X$ be a compact and connected submanifold.
We are interested in the field $\mathbb C(X,Y)$ of germs of meromorphic functions along $Y$, \ie equivalence classes of meromorphic functions defined on arbitrarily small open subsets of $X$ containing $Y$.

When the normal bundle of $Y$ in $X$ is ample, classical results by Andreotti and Hartshorne imply that the transcendence degree of $\mathbb C(X,Y)$  over $\mathbb C$ is bounded by the dimension of $X$. Recently, in \cite{MR4403217}, the first two of the authors of the present paper proved the existence of germs of smooth surfaces $S$
along a smooth compact curve $C$ isomorphic to $\mathbb P^1$ with arbitrary positive self-intersection without non-constant  meromorphic functions.

Our first main result is a simpler, yet more general, construction of germs of complex manifolds without non-constant meromorphic functions that generalizes the main result of \cite{MR4403217}.

\begin{THM}\label{THM:A}
    For any complex projective manifold $Y$ of dimension $n$ and any pair of natural numbers $(\ell, m)$ such that $m \ge 2n$ and $\ell \le m$, there exists a germ of complex manifold $X$ of dimension $m$ and containing $Y$  such that the normal bundle of $Y$ in $X$ is ample and the transcendence degree of $\mathbb C(X,Y)$ over $\mathbb C$ is equal to $\ell$. Moreover, when $\ell =0$ we can further guarantee that $X$ carries no irreducible subvarieties of dimension at least $n+1$ neither foliations/webs of arbitrary non-zero dimension/codimension.
\end{THM}

More informally, if a submanifold $Y\subset X$ with ample normal bundle has dimension smaller than its codimension then there are no restrictions on the transcendence degree
of $\C(X,Y)$ over $\C$ besides the upper bound $\dim X$.

In sharp contrast, when $Y$ is a smooth hypersurface with ample normal bundle on a complex manifold $X$ and $\dim Y \ge 2$, Rossi proved in \cite[Section 5, Theorem 3]{MR0176106} the
existence of a projective manifold $\overline X$ and an open subset $U \subset X$ containing $Y$ that can be identified with an open subset of $\overline X$.

\begin{question}\label{Q:wild guess}
    Let $X$ be a complex manifold and $Y \subset X$ a compact submanifold with ample normal bundle.
    Does $\dim X < 2 \dim Y$ (\ie $\codim Y< \dim Y)$  implies $\trdeg \mathbb C(X,Y) = \dim X$ ?
\end{question}

Theorem \ref{THM:A} and Question \ref{Q:wild guess} should be compared with a conjecture by Peternell \cite{MR3049294} which predicts that field $\C(X)$ of meromorphic functions on  a compact Kähler manifold $X$ admitting a subvariety $Z$ with ample normal bundle satisfies $\trdeg \mathbb C(X) \ge \dim Z + 1$. We do not know if we can construct examples of pairs $(X,Y)$ with $X$ compact Kähler with arbitrary $\trdeg \mathbb C(X,Y)$. Peternell's conjecture predicts that this is not the case.

One of the main ingredients of our proof of Theorem \ref{THM:A} is an old trick attributed to Haefliger and explored by Loeb-Nicolau, López de Medrano-Verjovsky, and Meersseman, see \cite[Introduction]{MR1760670} and references therein. In order to construct complex manifolds, it suffices to construct smooth differentiable manifolds transverse to
holomorphic foliations. Another key ingredients used in the proof of Theorem \ref{THM:A} are the very same results by Andreotti and Hartshorne mentioned above that served as motivation for this work.

Although we are not able to provide an answer to Question \ref{Q:wild guess}, our second main result shows that the method of proof of Theorem \ref{THM:A}
cannot be applied to produce counter-examples to it.

\begin{THM}\label{T:Van de Ven}
    Let $\F$ be a codimension $q$ foliation on $\mathbb P^n$ and let $Z \subset \mathbb P^n$ be a submanifold disjoint from $\sing(\F)$.
    If $q < 2 \dim Z$ and,  for every $z \in Z$, the intersection of the tangent spaces of $\F$ and $Z$ at $z$ has the expected
    dimension $\max \{ 0, \dim Z + \dim \F -n\}$ then $\F$ is algebraically integrable.
\end{THM}

As we will explain in Section \ref{S:Van de Ven}, Theorem \ref{T:Van de Ven} is strictly related to a classical result of Van de Ven \cite{MR0116361} characterizing linear subspaces of projective spaces. We also point out, that it is unclear what the optimal statement should be. It is conceivable that under the assumptions of Theorem \ref{T:Van de Ven} the foliation
$\F$ must defined by a linear projection (\ie $\F$ has degree zero). We actually prove this when $\F$ is a foliation by curves, see Theorem \ref{T:weakly transverse}.

\subsection{Acknowledgments} The authors thank Brazilian-French Network in Mathematics and CAPES/COFECUB Project Ma 932/19 “Feuilletages holomorphes et intéractions avec la géométrie”. 
Falla Luza acknowledges support from CNPq (Grant number 402936/2021-3).
Loray is supported by CNRS and  Centre Henri Lebesgue, program ANR-11-LABX-0020-0.
Pereira acknowledges support from CNPq (Grant number 301683/2019-0), and FAPERJ (Grant number E-26/200.550/2023). 
\section{Extension}

\subsection{Conventions}
Although we are mainly interested on submanifolds of complex manifolds, we will need to
consider the more general notion of complex analytic spaces in the sense of Grauert,
see for instance \cite[Appendix B]{MR0463157}. To wit,  subvarieties are closed complex analytic
subspaces without nilpotents in their structural sheaves, and submanifolds are smooth subvarieties.

\subsection{Infinitesimal neighborhoods}
Let $X$ be a complex manifold and let $Y\subset X$ be a closed and connected complex analytic subspace.
If $\mathcal I$ is the defining ideal sheaf of $Y$ then $n$-th infinitesimal neighborhood
of $Y$ is the complex analytic space $Y(n)$ with underlying topological space equal to $Y$ and structural sheaf equal to
\[
    \frac{\mathcal O_X}{\mathcal I^{n+1}} \, .
\]

There are natural inclusions $Y(n) \hookrightarrow Y(n+1)$ between infinitesimal neighborhoods of $Y$ induced by the natural
restriction morphisms
\[
    \frac{\mathcal O_X}{\mathcal I^{n+2}} \to \frac{\mathcal O_X}{\mathcal I^{n+1}} \, .
\]
The formal completion of $X$ along $Y$, denoted here by $Y(\infty)$, is the formal analytic space
obtained as the direct limit of $Y(n)$ when $n$ goes to $\infty$, \ie
\[
    Y(\infty) = \varinjlim Y(n) \, .
\]
In other words, $Y(\infty)$ is the ringed space with underlying topological space equal to $Y$ and with structural sheaf equal to the inverse limit
\[
    \mathcal O_{Y(\infty)} = \varprojlim \frac{\mathcal O_X}{\mathcal I^{n+1}} = \varprojlim \mathcal O_{Y(n)} \, .
\]

\subsection{The field of formal meromorphic functions}\label{SS:field}
Following \cite{MR251043}, we define the sheaf $\mathcal M_{Y(\infty)}$ of (formal) meromorphic functions on $Y(\infty)$ as the
sheaf associated to the presheaf  with sections over an open subset $U\subset Y$ equal to the total ring of fractions of $\mathcal O_{Y(\infty)}(U)$. We will denote the global sections of $\mathcal M_{Y(\infty)}$ by $\mathbb C(Y(\infty))$.

As we are assuming that $Y$ is connected $\mathbb C(Y(\infty))$ is a field, see \cite[Proposition 9.2]{MR2103516}.  Moreover,
we have  injective restriction morphisms
\begin{equation}\label{E:inclusions}
    \mathbb C(X) \hookrightarrow \mathbb C(X,Y) \hookrightarrow \mathbb C(Y(\infty)).
\end{equation}
When $Y$ is a point on a $n$-dimensional projective manifold $X$, then $\mathbb C(X,Y)$ is the quotient field of the ring
of convergent power series $\C\{x_1, \ldots, x_n\}$  while $\mathbb C(Y(\infty))$ is the quotient field of the ring of formal power series $\C[[x_1, \ldots, x_n]]$. Therefore, in general, both inclusions can have infinite transcendence degree.

In sharp contrast, Hironaka and Matsumura \cite[Theorem 3.3]{MR251043} proved the following result for subvarieties of projective spaces.

\begin{thm}\label{T:HironakaMatsumura}
    Let $Y$ be a closed and connected complex analytic subspace of $\mathbb P^n$, $n \ge 2$. If $\dim Y \ge 1$ then $\C(Y(\infty)) = \C(X,Y) = \C(\mathbb P^n)$.
\end{thm}

For an analytic version of Theorem \ref{T:HironakaMatsumura}, establishing the equality $\C(X,Y) = \C(\mathbb P^n)$ through Hartog's theorem and without the use of formal geometry,  see \cite[Theorem 3]{MR244516}.

When $\dim Y>0$ and the normal bundle of $Y$ is ample,  we can rephrase a result by Hartshorne \cite[Theorem 6.7 and Corollary 6.8]{MR232780}
to our setting as follows.

\begin{thm}\label{T:Hartshorne}
    Let $Y$ be a connected and compact complex analytic subspace of a complex manifold $X$. If $Y$ is locally a complete intersection of dimension at least one and the normal bundle $N_{Y/X}$ is ample then the transcendence degree over $\C$ of $\C(Y(\infty))$ is bounded by the dimension of $X$. Moreover, if $\trdeg \C(Y(\infty)) = \dim X$ then $\C(Y(\infty))$ is a finitely generated extension of $\C$.
\end{thm}

Although Theorem \ref{T:Hartshorne} was originally stated in the algebraic category, Hartshorne's proof works, as it is, in the context of complex manifolds considered here.

When the ambient manifold $X$ is projective (or more generally Moishezon) and $Y$ has ample normal bundle, Theorem \ref{T:Hartshorne} implies that $\C(Y(\infty))$ is a finite algebraic extension of $\mathbb C(X)$. Consequently, in this case, both inclusions in (\ref{E:inclusions}) are finite algebraic extensions. This fact combined with the proposition below,  due to Badescu and Schneider (see  \cite[Proposition 3.5]{MR1954055} or \cite[Proposition 10.17]{MR2103516}), imply the convergence of formal meromorphic functions defined on formal neighborhoods of subvarieties with ample normal bundle on projective manifolds.

\begin{prop}\label{P:Badescu Schneider}
    Let $X$ be a projective manifold  and $Y \subset X$ a connected subvariety. Then the algebraic closure of $\C(X)$ inside of $\C(Y(\infty))$
    is contained in $\mathbb C(X,Y)$.
\end{prop}

\begin{cor}\label{C:converge}
    Let $X$ be a projective manifold  and $Y \subset X$ a connected subvariety with ample normal bundle. Then   $\C(X,Y)=\C(Y(\infty))$.
\end{cor}

\begin{question}\label{Q:converge?}
    If the ambient space is not projective, just a small euclidean neighborhood of $Y$, is it true that  the ampleness of $N_Y$ implies
    the convergence of the formal meromorphic functions ?
\end{question}

Rossi's theorem (\cite{MR0176106}) gives a positive answer when $Y$ is a hypersurface of dimension at least two.

\subsection{Continuation of analytic subvarieties} The next result was proved by Rossi in \cite[Theorem 3.2]{MR244516} under the additional assumption $\dim (V \cap Y) + n =  \dim Y + \dim V$. The improvement below is due to Chow \cite[Corollary of Theorem 4]{MR257074}.  Unaware of Chow's improvement on Rossi's result, one of the authors of the present paper obtained an alternative proof of Rossi-Chow Theorem   in \cite[Theorem C]{MR3999055}.

\begin{thm}\label{T:Rossi}
    Let $Y \subset \mathbb P^n$ be an irreducible subvariety and  let $U\subset \mathbb P^n$ be a Euclidean neighborhood of $Y$.
    If $V\subset U$ is an irreducible subvariety of $U$ with  $\dim V + \dim Y > n$ and non-empty intersection with $Y$ then there exists a projective subvariety $\overline V \subset \mathbb P^n$ such that $V$ is an irreducible component of $\overline V \cap U$.
\end{thm}

Bost \cite[Theorem 3.5]{MR1863738}, Bogomolov-McQuillan \cite{MR3644242},  Campana-Paun \cite[Theorem 1.1]{MR3949026}, Druel \cite[Proposition 8.4]{MR3742759}, and others, used variants of Hartshorne's argument to establish algebraicity criteria for leaves of foliations on projective varieties. See also the work of Chen \cite[Theorem A]{MR2957623} for an algebraicity criterion for formal subschemes of projective schemes. Here we point out that we can deduce from Theorem  \ref{T:Hartshorne} a variant of Theorem \ref{T:Rossi}.

\begin{prop}\label{P:Rossi variant}
    Let $X$ be a projective manifold, let $Y \subset X$ be an irreducible locally complete intersection analytic subspace/subscheme with ample normal bundle, and let $U \subset X$ be a small Euclidean neighborhood of $Y$. If $V \subset U$ is a closed, irreducible, smooth subvariety of $U$ such that
    \begin{enumerate}
        \item\label{I:Rossi variant 1} $\dim V + \dim Y > \dim X$ and
        \item\label{I:Rossi variant 2} $\dim (V \cap Y) + \dim X =  \dim Y + \dim V$
    \end{enumerate}
    then there exists a projective subvariety $\overline V \subset X$ such that $V$ is an irreducible component of $\overline V \cap U$.
\end{prop}
\begin{proof}
    Condition (\ref{I:Rossi variant 1}) implies that $\dim V \cap Y >0$.
    Condition (\ref{I:Rossi variant 2}) implies that the codimension of $V\cap Y$ in $V$ equals the codimension of $Y$ in $X$.
    Therefore $V\cap Y$ is a locally complete intersection analytic subspace of $V$ with ample normal bundle and positive dimension. Consider the Zariski closure $\overline V \subset X$ of $V$ and notice that 
    $\trdeg \C(\overline V)=\trdeg \C(V\cap Y(\infty))$. On the other hand, we can apply Theorem \ref{T:Hartshorne}
    to $V\cap Y\subset V$ and get that $\trdeg \C(V\cap Y(\infty))\le \dim(V)$, and conclude.
\end{proof}

\section{Construction}

We carry out the proof of Theorem \ref{THM:A} in this section.

\subsection{Weakly transverse submanifolds}
Let $\F$ be a singular holomorphic foliation on a complex manifold $W$. We will say that
a submanifold $Z\subset W$ is weakly transverse to  $\F$ if
\begin{enumerate}
    \item the singular set of $\F$ does not intersect $Z$; and
    \item for every point $z \in Z$, the tangent space of $Z$ at $z$ intersects the tangent space of $\F$ at $z$
    only at $0$.
\end{enumerate}

The relevance of the concept to our discussion is put in evidence by our next result.

\begin{prop}\label{P:leaf space}
    Let $\F$ be a singular holomorphic foliation on a complex manifold $W$.
    If $Z \subset W$ is a compact submanifold weakly transverse to $\F$ then there exist
    a neighborhood $U$ of $Z$; a complex manifold $X$ of dimension equal to $\dim W- \dim \F$; and
    a holomorphic submersion $\pi : U \to X$ such that
    \begin{enumerate}
        \item the leaves of  $\restr{\F}{U}$ coincide with the fibers of $\pi$; and
        \item the morphism $\pi$ maps  $Z$ isomorphically to a submanifold $Y$ of $X$  with normal bundle isomorphic to
        the quotient of $N_{Z/W}$ by the image of $\restr{T_{\F}}{Z}$ inside it; and
        \item the field of germs of meromorphic functions $\mathbb C(X,Y)$ is mapped by $\pi^*$ isomorphically
        onto the field of germs of meromorphic first integrals of the germ of foliation $\restr{\F}{(W,Z)}$.
    \end{enumerate}
\end{prop}
\begin{proof}
    Let $U$  be equal to a sufficiently small tubular neighborhood of $Z$ in $W$.  The weakly transversality between $\F$ and $Z$ implies that the leaves of $\restr{\F}{U}$ are closed subvarieties of $U$. Therefore the leaf space $X=U/\F$ of $\restr{\F}{U}$ is a Hausdorff complex manifold. The quotient morphism $\pi: U \to X$ has the sought properties.
\end{proof}

\subsection{Formal meromorphic functions versus rational first integrals} The field of meromorphic functions
of the leaf space constructed in Proposition \ref{P:leaf space} is described by the lemma below.

\begin{lemma}\label{L:field of functions}
    Notation and the assumptions as in Proposition \ref{P:leaf space}. Further assume that $W = \mathbb P^m$. If $Y(\infty) \subset X$ is the formal completion of $Y= \pi(Z)$ in $X$ then $\mathbb C(Y(\infty))$ is isomorphic to the field of rational first integrals of the foliation $\F$.
\end{lemma}
\begin{proof}
    Direct consequence of Hironaka-Matsumura result discussed in Subsection \ref{SS:field}. Indeed, the pull-back of $\mathbb C(Y(\infty))$
    under $\pi$ is the subfield of $\mathbb C(Z(\infty))$ formed by the formal meromorphic first integrals of $\restr{\F}{Z(\infty)}$. By Hartshorne result, $\mathbb C(Z(\infty))$ coincides with the field  $\mathbb C(W) = \mathbb C(\mathbb P^m)$ the field of rational functions on $W$. Consequently, we get an isomorphism between $\mathbb C(Y(\infty))$ and the field of rational first integrals of $\F$.
\end{proof}

\begin{remark}\label{R:super}
    Likewise, when $\F$ is a global foliation on $\mathbb P^n$,  foliations/webs living in the leaf space of $\restr{\F}{U}$ correspond to  global foliations/global webs invariant by $\F$. This again is a consequence of Hartshorne's result.
\end{remark}

The situation for subvarieties of the leaf space is slightly different, and it is not true that their pre-iamges under the quotient morphism  can always be globalized.

\begin{lemma}\label{L:subvarieties}
    Notation and the assumptions as in Proposition \ref{P:leaf space}. Further assume that $W = \mathbb P^m$. If $S \subset X = U/\F$ is a subvariety intersecting $Y= \pi(Z)$ such that $\dim S + \dim Y > \dim X$ then $\pi^{-1}(S)$ is contained in a $\F$-invariant subvariety $\overline S$ of $\mathbb P^m$ with $\dim \overline S = \dim S + \dim \F$.
\end{lemma}
\begin{proof}
    Direct consequence of  Theorem \ref{T:Rossi}.
\end{proof}

\subsection{Existence of weakly transverse submanifolds} As shown below, the existence of weakly transverse submanifolds, under suitable numerical assumptions, 
is easy consequence of Kleiman's transversality of a general translate.

\begin{lemma}\label{L:existence}
    Let $\F$ be a foliation by curves on $\mathbb P^{m+1}$ with isolated singularities.
    Let $Y \subset \mathbb P^{m+1}$ be a projective submanifold of dimension $n$.
    If $m\ge 2n$ then $Y$ is weakly transverse to $g^* \F$ for any general $g \in \Aut(\mathbb P^{m+1})$.
\end{lemma}
\begin{proof}
    The proof is a simple application of Kleiman's transversality of a general translate, \cite[Theorem 2]{MR360616}.
    Indeed, we can identify the projectivization of the tangent bundle of $Y$, $\mathbb P(T_Y)$ (space of lines on $T_Y$),
    with a submanifold $\tau(Y)$ of $\mathbb P T_{\mathbb P^{m+1}}$ of dimension $2n-1$. Likewise, we consider $\tau(\F) \subset \mathbb P(T_{\mathbb P^m+1})$ as the Zariski closure of the tangent lines of $\restr{\F}{\mathbb P^{m+1} - \sing(\F)}$. As such $\tau(\F) \subset \mathbb P (T_{\mathbb P^{m+1}})$ is a subvariety of dimension $m+1$. By assumption
    \[
        \dim \tau(Y) + \dim \tau(\F) < \dim \mathbb P(T_{\mathbb P^n}) \, .
    \]
    Since the natural action of $\Aut(\mathbb P^{m+1})$ on $\mathbb P (T_{\mathbb P^{m+1}})$ is transitive, we can apply
    \cite[Theorem 2]{MR360616} to guarantee that, for a general $g \in \Aut(\mathbb P^{m+1})$, $g^*\tau(\F)=\tau(g^*\F)$ is disjoint from $\tau(y)$. The lemma follows.
\end{proof}

\subsection{Field of rational first integrals of foliations on projective manifolds}

Let $\F$ be a foliation on a projective manifold $X$. According to \cite{MR2223484}, there exists a unique foliation by algebraic leaves $\G$ containing $\F$ and such that $\mathbb C(\F) = \mathbb C(\G)$.

\begin{ex}\label{Ex:linear}
    Let $\lambda_1, \ldots, \lambda_{m+1}$ be complex numbers and consider the following vector field on $\mathbb C^{m+1}$
    \[
        v = \sum_{i=1}^{m+1} \lambda_i x_i \frac{\partial}{\partial x_i} \, .
    \]
    The foliation $\F$ on $\mathbb P^{m+1}$ defined by $v$ has Zariski closure $\G$ of dimension equal to
    the dimension of $\mathbb Q$-vector subspace of $\mathbb C$ generated by $\lambda_1, \ldots ,\lambda_{m+1}$.
    In particular, choosing appropriately the complex numbers $\lambda$ we have one dimensional foliations on $\mathbb P^{m+1}$
    having field of rational first integrals of any transcendence degree between $0$ and $m$.
\end{ex}

\subsection{Proof of Theorem \ref{THM:A}}
Let $\F$ be a one-dimensional foliation $\mathbb P^{m+1}$ with field of rational first integrals of transcendence degree $\ell$ over $\C$. For instance, we can take a foliation $\F$ as in Example \ref{Ex:linear}.
Embed $Y$ in $\mathbb P^{m+1}$. According to Lemma \ref{L:existence} we can assume that $Y$ is weakly transverse to $\F$. Proposition \ref{P:leaf space} implies the existence of an Euclidean open subset $U \subset \mathbb P^{m+1}$ containing $Y$ such that the leaf space $X = U / \restr{\F}{U}$
has $\C(X)$ equal to $\C(\restr{\F}{U})$. Lemma \ref{L:field of functions} implies that $\C(X)$ equals to $\C(\F)$. This shows that the leaf space $X = U / \restr{\F}{U}$ is a complex manifold of dimension $m$ containing $Y$ and with $\trdeg \mathbb C(X,Y) = \ell$ as claimed.

The last claim concerning the existence of manifolds with $\ell =0$ follows from Remark \ref{R:super} and Lemma \ref{L:subvarieties} combined with \cite[Theorem 1]{MR2862041} which guarantees that the very general one-dimensional foliation on $\mathbb P^{m+1}$ of degree at least two is not tangent to any other foliation or web and does not leave invariant any algebraic subvariety.
\qed

\section{Algebraicity}\label{S:Van de Ven}

In this section we will present a proof of Theorem \ref{T:Van de Ven} from the Introduction.
Actually Theorem \ref{T:Van de Ven} is the combination of Theorems \ref{T:transverse} and \ref{T:weakly transverse}
below.

\subsection{Submanifolds transverse to foliations} We start by treating submanifolds transverse to foliations. 

\begin{thm}\label{T:transverse}
    Let $\F$ be a codimension $q$ foliation on a projective space $\mathbb P^n$.  If $Z \subset \mathbb P^n$ is a submanifold
    disjoint from $\sing(\F)$ and transverse to $\F$ then $\dim Z = q$,  $Z$ is a linear subspace, and the foliation
    $\F$ has degree zero.
\end{thm}
\begin{proof}
    Let $i: Z \to \mathbb P^n$ be the natural inclusion and let $\omega \in H^0(\mathbb P^n, \Omega^q_{\mathbb P^n} \otimes \det N_{\F})$  be a $q$-form
    defining $\F$. Observe that $\dim Z \ge q = \codim \F$ by the definition of transversality.
    Moreover, $i^* \omega \in H^0(Z, \Omega^q_{Z} \otimes \restr{\det N_{\F}}{Z})$ is an everywhere non-zero
    twisted $q$-form on $Z$. If $\dim Z \ge q+1$ then $i^*\omega$ defines a smooth foliation on $Z$ of positive dimension. Bott's vanishing
    theorem implies that $(\restr{\det N_{\F}}{Z})^{q+1} = 0$. But $\restr{\det N_{\F}}{Z}$ is an ample line-bundle and  $q+1\le \dim Z$, hence
    $(\restr{\det N_{\F}}{Z})^{q+1} \neq 0$. We obtain in this way  a contradiction that establishes that $\dim Z=q$.
    Therefore the non-vanishing of $i^*\omega$ implies that the canonical sheaf $\Omega^q_Z$ is isomorphic to $\restr{\det N^*_{\F}}{Z} = \mathcal O_Z(- \deg(\F) - q - 1)$.
    Kobayashi-Ochiai Theorem \cite[Corollary of Theorem 1.1]{MR316745} implies that $Z = \mathbb P^{q}$ , $\deg(\F) =0$ and $\deg (Z) = 1$.
\end{proof}

\begin{remark}
    Let $\F$ be a codimension $q$ foliation on a projective manifold $X$. If $Z \subset X$ is a submanifold with ample normal bundle, disjoint from $\sing(\F)$, and transverse to $\F$ then proof of Theorem \ref{T:transverse} shows that $\dim Z=q$. Hence $\restr{T_{\F}}{Z} \simeq N_Z$ is ample and we can apply  \cite[Theorem 3.5]{MR1863738} or \cite{MR3644242} to guarantee that every leaf of $\F$ passing through the points of $Z$ is algebraic. Consequently, by dimension reasons, every leaf of $\F$ is algebraic and $\F$ is algebraically integrable.
\end{remark}

We present below a proof of a classical result by Van de Ven.

\begin{cor}
    Let $X$ be a submanifold of $\mathbb P^n$. If the normal sequence
    \[
        0 \to T_X \to \restr{T_{\mathbb P^n}}{X} \to N_X \to 0
    \]
    splits then $X$ is a linear subspace of $\mathbb P^n$.
\end{cor}
\begin{proof}
    Let $q$ be the codimension of $X$.
    The existence of a splitting $\varphi: N_X \to \restr{T_{\mathbb P^n}}{X}$ gives us for every $x \in X$, a unique linear $\mathbb P^q$
    through $x$ with tangent space at $x$ equal to $\varphi(N_X(x))$. If $U$ is a sufficiently nieghborhood of $X$ then this family of $\mathbb P^q$'s will not intersect at $U$, and will thus define a codimension $q$ foliation $\F_U$ on $U$. Theorem \ref{T:Hartshorne} implies that $\F_U$ extends to a foliation $\F$ on $\mathbb P^n$. By construction $\F$ is transverse to $X$. Theorem \ref{T:transverse} implies the result.
\end{proof}

\subsection{Submanifolds weakly trasnverse to foliations} Our last result treats the case of submanifolds weakly transverse to foliations. 

\begin{thm}\label{T:weakly transverse}
    Let $\F$ be a codimension $q$ foliation on $\mathbb P^n$. If $Z$ is a submanifold
    weakly transverse to $\F$ and $q <2\dim Z$ then $\F$ is algebraically integrable. Moreover, if $\dim \F =1$ then $\deg(\F)=0$.
\end{thm}
\begin{proof}
    Let $U \subset \mathbb P^n$ be a sufficiently small Euclidean neighborhood of $Z$.
    Let $\pi : U \to U/\F$ be the quotient of $U$ by $\restr{\F}{U}$.  Observe that $V=\pi^{-1}(\pi(Z))$ is a
    closed subvariety of $U$ of dimension
    \begin{equation} \label{E:dimensao V}
        \dim \F + \dim Z.
    \end{equation}
    Hence our assumptions imply that $\dim Z + \dim V = 2 \dim Z + \dim \F > n$. Theorem \ref{T:Rossi} implies that $\overline{V}$, the Zariski closure of $V$, and $V$ have  the same dimension. Since $V$ is invariant by $\restr{\F}{U}$, $\overline{V}$ is also invariant by $\F$.

    Let $\G$ be the unique foliation by algebraic leaves containing $\F$ and with $\mathbb C(\F) = \mathbb C(\G)$ given by \cite{MR2223484}.
    After replacing $Z$ by $gZ$ for a sufficiently general automorphism $g \in \Aut(\mathbb P^n)$, we can assume that the Zariski closure of a  leaf of $\F$ through a general point $z$ of $Z$ coincides  with a leaf of $\G$. Therefore
    \begin{equation}\label{E:dimensao V barra}
        \dim \overline{V} = \dim \G + \dim Z.
    \end{equation}
    Since $\dim V = \dim \overline{V}$, we can combine  Equations (\ref{E:dimensao V}) and (\ref{E:dimensao V barra}) to establish the algebraicity of the leaves of $\F$, in other terms $\F = \G$.

    In the one dimensional case, weakly transversality implies the existence of a nowhere zero section of $N_Z\otimes \omega_{\F}|_Z$, thus $c_{\mathrm{cod}(Z)}(N_Z\otimes \omega_{\F})=0$.

    Since $N_Z$ is an ample vector bundle, and the Chern classes of ample vector bundles are strictly positive according to \cite[Theorem 2.5]{MR297773}, we deduce that
    $\omega_{\F}$ is not nef. Since $\omega_{\F}|_Z=\mathcal{O}_Z(\deg(\F) - 1)$ we necessarily have $\deg(\F)=0$ as claimed.
\end{proof}

\begin{remark}
    We believe that a similar statement should hold true if $\mathbb P^n$ is replaced by an arbitrary projective manifold and $Z$ is replaced by a submanifold with ample normal bundle.  
\end{remark}

\bibliography{references}{}
\bibliographystyle{amsplain}

\end{document}